\documentstyle[amstex,amssymb,A4,11pt]{article}
\makeatletter

\def\@@N{{\mathrm I \mkern-2.5mu \nonscript\mkern-.5mu \mathrm I
\mkern-5.5mu\mathrm  N}}
\def\@@Z{{
  \setbox0\hbox{\sf Z}
  \setbox1\hbox{\mathrm\kern.05\wd0
                \rlap{\vrule height.93\ht0 depth-.75\ht0 width.056\wd0 }%
                \kern-.13\wd0 \copy0 \kern-.6\wd0 \copy0 \kern-.1\wd0
                \llap{\vrule height.25\ht0 depth\z@ width.056\wd0}%
                \kern.05\wd0}
  \mathchoice{\copy1}{\copy1}{\mit Z\mkern-8mu Z}{\mit Z\mkern-7.5mu Z} }}
\def\@@insvline#1#2{{\setbox0\hbox{\m@th$#1\mathrm I$}
  \rlap{\m@th$#1 \mkern 5mu
  \vrule height.92\ht0 depth-.05\ht0 width.09\ht0 $}
  {\mathrm #2} }}
\def\@@Q{\mathpalette\@@insvline{Q}}
\def\@@R{{\mathrm I \mkern-2.5mu \nonscript\mkern-.5mu \mathrm R}}
\def\@@C{\mathpalette\@@insvline{C}}
\def\@@K{{\mathrm I \mkern-2.5mu \nonscript\mkern-.5mu\mathrm  K}}
\def\@@D{{\mathrm I \mkern-2.5mu \nonscript\mkern-.5mu\mathrm D}}

\def\callig#1{
  \ifcat#1a
    \ifnum`#1=\uccode`#1 {\cal #1}
    \else                #1
    \fi
  \else                  #1
  \fi}

\def\varbbbold#1{
  \ifcat#1a
    \ifnum`#1=\uccode`#1 \csname @@#1\endcsname
    \else                #1
    \fi
  \else                  #1
  \fi}

\def\eulerfraktur#1{
  \ifcat#1a {\frak #1}
  \else                  #1
  \fi}

\newif\if@@excloccurred
\newif\if@@questoccurred

\let\@@quest=? \catcode`?=\active
\let\@@excl=!  \catcode`!=\active \let\fact\@@excl

\def!{\protect\p@@doexcl}
\def?{\protect\p@@doquest}

\def\p@@doexcl{\let\@@xeq\relax\ifmmode\def\@@xeq{\futurelet\@@next\@@doexcl}
                               \else\@@excl\fi\@@xeq}
\def\p@@doquest{\let\@@xeq\relax\ifmmode\def\@@xeq{\futurelet\@@next\@@doquest}
                                \else\@@quest\fi\@@xeq}
\def\@@doexcl{\let\@@x@q\relax
     \if@@excloccurred \@@excloccurredfalse
                       \ifx\@@next\@sptoken\@@excl\@@excl
                       \else\let\@@x@q\doexclexcl\fi
\else\if@@questoccurred\@@questoccurredfalse
                       \ifx\@@next\@sptoken\@@quest\@@excl
                       \else\let\@@x@q\doquestexcl\fi
\else\ifx\@@next\@sptoken\@@excl
\else\ifx\@@next!\@@excloccurredtrue
\else\ifx\@@next?\@@excloccurredtrue
\else\let\@@x@q\doexcl\fi\fi\fi\fi\fi\@@x@q}

\def\@@doquest{\let\@@x@q\relax
     \if@@excloccurred \@@excloccurredfalse
                       \ifx\@@next\@sptoken\@@excl\@@quest
                       \else\let\@@x@q\doexclquest\fi
\else\if@@questoccurred\@@questoccurredfalse
                       \ifx\@@next\@sptoken\@@quest\@@quest
                       \else\let\@@x@q\doquestquest\fi
\else\ifx\@@next\@sptoken\@@quest
\else\ifx\@@next!\@@questoccurredtrue
\else\ifx\@@next?\@@questoccurredtrue
\else\let\@@x@q\doquest\fi\fi\fi\fi\fi\@@x@q}

\def\doexcl{\bold}        \def\doquest{\callig}
\def\doexclexcl{\varbbbold} \def\doexclquest{\boldsymbol}
\def\doquestexcl{\Bbb}    \def\doquestquest{\eulerfraktur}
\makeatother

\newcommand {\art}[6]{{\sc #1:} {#2.} {\em #3} {\bf #4}{(#5),} {#6.}}
\newcommand {\artt}[6]{{\sc #1:} {#2} {\em #3} {\bf #4}{(#5),} {#6.}}
\newcommand {\bo}[5]{{\sc #1:} {``#2."} {#3,} {#4} {#5.}}
\newcommand {\samp}[8]{{\sc #1:} {#2,} {{\em in} ``#3,"} {(#4, Ed.),}
                                  {pp. #5,} {#6,} {#7} {#8.}}

\newcommand {\unfart}[3]{{\sc #1:} {#2,} {(#3).}}

\newcommand \force {{\hspace{0.4mm}}{\rule{0.1mm}{2.4mm}}
                       {\hspace{0.5mm}}{\rule{0.1mm}{2.4mm}}
               {\rule [1.2mm]{2.3mm}{0.1mm}}{\hspace{0.4mm}}}
\newcommand \reals {[\omega ]^\omega }

\newcommand \Null {{\hspace{-0.7mm}}{{}_\circ}}
\newcommand \Nullem {{}_\circ}

\newcommand \pp {{X}}
\newcommand \fD {\frak{D}}

\newcommand \Padname {{\operatorname{P}}_{\tilde{\frak{D}}}}

\newcommand \oQ {\operatorname{\bold{Q}}}

\newcommand \add {\operatorname{\bold{add}}}
\newcommand \cov {\operatorname{\bold{cov}}}

\newcommand \fW {\frak{W}}

\newcommand \fc {\frak{c}}
\newcommand \fI {\frak{i}}

\newcommand \fJ {\frak{J}}
\newcommand \ft {\frak{t}}
\newcommand \fh {\frak{h}}
\newcommand \fs {\frak{s}}
\newcommand \fS {\frak{S}}

\newcommand \fr {\frak{r}}
\newcommand \fd {\frak{d}}
\newcommand \fb {\frak{b}}
\newcommand \fR {\frak{R}}
\newcommand \fO {\frak{O}}
\newcommand \fu {\frak{u}}
\newcommand \fo {\frak{o}}
\newcommand \fp {\frak{p}}
\newcommand \di {\operatorname{\bold{d}}}
\newcommand \fH {\frak{H}}

\newcommand \J  {{\cal{J}}}

\newcommand \ceq {\sqsubseteq}
\newcommand \puni {\sqcap}
\newcommand \pint {\sqcup}
\newcommand {\open}[2]{{(#1,} {#2)^{\omega}}}

\newcommand \parto {(\omega)^\omega}
\newcommand \partou {(\omega)^{\underline{\omega}}}
\newcommand \partf {(\omega)^{<{\omega}}}
\newcommand \partX {(X)^\omega}

\newcommand \NN {(!!N)}
\newcommand \mdom {\operatorname{dom}}
\newcommand \mmin {\operatorname{min}}
\newcommand \mmax {\operatorname{max}}
\newcommand \mMin {\operatorname{Min}}

\newcommand \ass {\stackrel{*}{=}}
\newcommand \seg {\preceq}

\newcommand \MA {\operatorname{\bold{MA}}}

\newcommand \mao {{\em{mao}}}
\newcommand \oto {{\omega}^{\omega}}
\newcommand \Lort {\angle}
\newcommand \ort {\bot}
\newcommand \comp {\parallel}
\newtheorem {nummer}{ }[section]
\newtheorem {thm}[nummer]{T\footnotesize{HEOREM}}
\newtheorem {lm}[nummer]{L\footnotesize{EMMA}}

\newtheorem {fct}[nummer]{F\footnotesize{ACT}}

\newtheorem {cor}[nummer]{C\footnotesize{OROLLARY}}

\newcommand \rmkone {{\bf{R{\footnotesize{EMARK 1}}:\hspace*{3mm}}}}
\newcommand \rmktwo {{\bf{R{\footnotesize{EMARK 2}}:\hspace*{3mm}}}}

\newcommand \proof {{\bf{P{\footnotesize{ROOF}}:\hspace*{3mm}}}}

  \def\eop{{\unskip\nobreak\hfil\penalty50\hskip8mm\hbox{}
  \nobreak\hfil
  {$\boldsymbol{\dashv}$}\parfillskip=0mm \par\smallskip}}
\newcommand \eof {${\ }_\dashv$\hfill}

\begin{document}
\begin{center}
      {\Large{\sc{On shattering, splitting and reaping partitions}}}
\end{center}
\smallskip

\begin{center}
\small{Lorenz Halbeisen}\footnote{The author
wishes to thank the {\em Swiss National Science Foundation\/} for
supporting him.}\\
\small{Universit\'e de Caen}\\
\small{France}
\end{center}
\medskip

{\small {\bf Keywords:} Cardinal invariants, partition properties,
dual-Mathias forcing.}\hfill\smallskip

{\small {\bf MS-Classification:} 03E05, 03E35, 03C25,
04A20, 05A18.}\hfill\bigskip

\begin{abstract}
In this article we investigate the dual-shattering cardinal $\fH$,
the dual-splitting cardinal $\fS$ and the dual-reaping cardinal
$\fR$, which are dualizations of the
well-known cardinals $\fh$ (the shattering cardinal, also known as
the distributivity number of ${\cal P}(\omega )/\mbox{{\em{fin}}}$),
$\fs$ (the splitting number) and $\fr$ (the reaping number).
Using some properties of the ideal $\fJ$ of
nowhere dual-Ramsey sets, which is
an ideal over the set of partitions of
$\omega$, we show that $\add (\fJ)=
\cov (\fJ)=\fH$. With this result we can
show that $\fH >\omega_1$ is consistent with ZFC and as a corollary
we get the relative consistency of $\fH >\ft$, where
$\ft$ is the tower number. Concerning $\fS$ we show
that $\cov ({\cal M})\leq \fS$
(where ${\cal M}$ is the ideal of the meager sets).
For the dual-reaping cardinal $\fR$ we get $\fp
\leq\fR\leq\fr$ (where $\fp$ is the pseudo-intersection number) and
for a modified dual-reaping number $\fR'$ we get $\fR'\leq\fd$
(where $\fd$ is the dominating number). As a consistency result we
get $\fR <\cov ({\cal M})$.
\end{abstract}

\bigskip

\section{The set of partitions}

A {\em partial partition\/} $\pp$ (of $\omega$) consisting of pairwise
disjoint, nonempty sets, such that $\mdom(\pp):=\bigcup \pp
\subseteq \omega$.
The elements of a partial partition
$\pp$ are called the blocks of $\pp$ and $\mMin(\pp)$ denotes
the set of the least elements of the blocks of $\pp$.
If $\mdom(\pp)=\omega$, then $\pp$ is called a {\em partition}.
$\{\omega\}$ is the partition such that each block is a singleton
and $\{\{\omega\}\}$ is the partition containing only one block.
The set of all partitions containing infinitely (resp.\,finitely) many
blocks is denoted by $\parto$ (resp.\,$\partf$).
By $\partou$ we denote the set of all infinite partitions
such that at least one block is infinite.
The set of all partial partitions with $\mdom(\pp)\in
\omega$ is denoted by $\NN$.
\hfill\smallskip

Let $\pp_1,\pp_2$ be two partial partitions. We say that $\pp_1$ is
{\em coarser\/} than $\pp_2$, or that
$\pp_2$ is {\em finer\/} than $\pp_1$, and write $\pp_1\ceq\pp_2$ if
for all blocks $b\in\pp_1$ the set $b\cap\mdom(\pp_2)$ is the union
of some sets $b_i\cap\mdom(\pp_1)$, where each $b_i$ is a
block of $\pp_2$. (Note that if $\pp_1$ is coarser than $\pp_2$,
then $\pp_1$ is in a natural way also contained in $\pp_2$.)
Let $\pp_1\puni\pp_2$ denotes the finest
partial partition which is coarser than $\pp_1$ and $\pp_2$
such that $\mdom(\pp_1\puni\pp_2)=
\mdom(\pp_1)\cup\mdom(\pp_2)$. Similarly
$\pp_1\pint\pp_2$ denotes the coarsest partial partition which is finer
than $\pp_1$ and $\pp_2$ such that
$\mdom(\pp_1\pint\pp_2)=\mdom(\pp_1)\cup\mdom(\pp_2)$.
\hfill\smallskip

If $f$ is a finite subset of $\omega$, then
$\{f\}$ is a partial partition with $\mdom(\{f\})=f$. For two partial
partitions $\pp_1$ and $\pp_2$ we write $\pp_1\ceq^*\pp_2$ if there is a
finite set $f\subseteq\mdom(\pp_1)$ such that $\pp_1\puni\{f\}\ceq\pp_2$
and say that $\pp_1$ is coarser${}^*$ than $\pp_2$.
If $\pp_1\ceq^*\pp_2$ and $\pp_2\ceq^*\pp_1$ then we write
$\pp_1\ass\pp_2$. If $X\ass\{\omega\}$, then $X$ is called {\em
trivial}.\hfill\smallskip

%
Let $\pp_1,\pp_2$ be two partial partitions. If each block of $\pp_1$
can be written as the intersection of a block of $\pp_2$ with
$\mdom(\pp_1)$, then we write $\pp_1\seg\pp_2$.
Note that $\pp_1\seg\pp_2$ implies $\mdom(\pp_1)\subseteq\mdom(\pp_2)$.
\hfill\smallskip

We define a topology on the set of partitions as follows. Let
$X\in\parto$ and $s\in\NN$ such that $s\ceq X$, then
${\open{s}{X}}:=\{Y\in\parto :s\seg Y\wedge Y\ceq X\}$ and $\partX
:={\open{\emptyset}{X}}$. Now let the basic open sets on $\parto$ be the
sets ${\open{s}{X}}$ (where $X$ and $s$ as above). These sets are called
the {\em dual Ellentuck neighborhoods.} The topology induced by the dual
Ellentuck neighborhoods is called the {\em dual Ellentuck topology\/}
(cf.~\cite{dual}).

\section{On the dual-shattering cardinal $\fH$}\label{sec:add}

\subsection*{Four cardinals}

We first give the definition of the dual-shattering cardinal $\fH$.
\hfill\smallskip

Two partitions $\pp_1,\pp_2\in\parto$ are called {\em{almost
orthogonal}\/} ($\pp_1\ort_* \pp_2$) if $\pp_1\puni\pp_2\not\in\parto$,
otherwise they are {\em{compatible}} ($\pp_1\comp\pp_2$). If $\pp_1
\puni\pp_2=\{\{\omega\}\}$, then they are called {\em orthogonal\/}
($\pp_1\ort\pp_2$). We say
that a family ${\cal A}\subseteq\parto$ is {\em{maximal almost orthogonal}\/}
(\mao) if ${\cal A}$ is a maximal family of pairwise almost orthogonal
partitions. A family ${\cal H}$ of {\mao} families of partitions
{\em{shatters}\/} a partition $\pp\in\parto$, if there are $H\in {\cal H}$ and
two distinct partitions in $H$ which are both compatible with $\pp$. A family of
{\mao} families of partitions is {\em{shattering}\/} if it shatters each member
of $\parto$. The dual-shattering cardinal $\fH$ is the least cardinal
number $\kappa$, for which there exists a shattering family of
cardinality $\kappa$.\hfill\smallskip

One can show that $\fH\leq\fh$ and $\fH\leq\fS$ (cf.\,\cite{dia}),
(where $\fS$ is the dual-splitting cardinal).

\subsubsection*{Two cardinals related to the ideal of nowhere dual-Ramsey
sets}

Let $C\subseteq\parto$ be a set of partitions, then we say that $C$ has
the {\em dual-Ramsey property\/} or that $C$ is {\em dual-Ramsey\/}, if
there is a partition $X\in\parto$ such that $\partX\subseteq C$ or
$\partX\cap C=\emptyset$. If the latter case holds, we also say that $C$
is {\em dual-Ramsey}$\Nullem$. If for each dual Ellentuck neighborhood
${\open{s}{Y}}$ there is an $X\in{\open{s}{Y}}$ such that
${\open{s}{X}}\subseteq C$ or ${\open{s}{X}}\cap C=\emptyset$, we call
$C$ {\em completely dual-Ramsey.\/} If for each dual Ellentuck
neighborhood the latter case holds, we say that $C$
is {\em nowhere dual-Ramsey}.\\
\rmkone
In \cite{dual} it is proved,
that a set is completely dual-Ramsey if and only if it has the Baire
property and it is
nowhere dual-Ramsey if and only if it is meager
with respect to the dual Ellentuck topology. From this it follows, that
a set is nowhere dual-Ramsey if and only if the complement contains
a dense and open subset (with respect to the dual Ellentuck topology).
\hfill\smallskip

Let $\fJ$ be set of partitions which are completely dual-Ramsey$\Null$.
The set $\fJ \subseteq {\cal P}(\parto)$ is an ideal which is not prime.
The cardinals $\add(\fJ)$ and $\cov(\fJ)$ are two cardinals related
to this ideal.\\
$\add (\fJ)$ is the smallest cardinal $\kappa$ such that there
exists a family ${\cal F}=\{J_\alpha\in\fJ :\alpha <\kappa\}$ with
$\bigcup{\cal F}\not\in\fJ$.\\
$\cov (\fJ)$ is the smallest cardinal $\kappa$ such that there
exists a family ${\cal F}=\{J_\alpha\in\fJ :\alpha <\kappa\}$ with
$\bigcup{\cal F}=\parto$.\\
Because $\parto\not\in\fJ$, it is clear that
$\add(\fJ)\leq\cov(\fJ)$. Further it is easy to see that
$\omega_1\leq\add(\fJ)$. In the next section we will show that
$\add(\fJ)=\cov(\fJ)$.

\subsubsection*{The distributivity number $\di (\fW)$}

A complete Boolean algebra $\langle B,\leq\rangle$ is called
$\kappa$-distributive, where $\kappa$ is a cardinal, if and only if for
every family $\langle u_{\alpha i}:i\in I_\alpha,\alpha<\kappa\rangle$ of
members of $B$ the following holds:
$$\prod\limits_{\alpha <\kappa}\sum\limits_{i\in I_\alpha}u_{\alpha i}=
\sum\limits_{f\in\prod\limits_{\alpha<\kappa}I_\alpha}
\prod\limits_{\alpha<\kappa}u_{\alpha
f(\alpha)}.$$
It is well known (cf.\,\cite{Jech}) that for a forcing notion $\langle
P,\leq\rangle$ the following statements are equivalent:
\begin{itemize}
\item r.o.($P$) is $\kappa$-distributive.
\item The intersection of $\kappa$ open dense sets in $P$ is dense.
\item Every family of $\kappa$ maximal anti-chains of $P$ has a common
refinement.
\item Forcing with $P$ does not add a new subset of $\kappa$.
\end{itemize}

Let $\J$
be the ideal of all finite sets of $\omega$ and
let $\langle {\parto}/{\J},\leq\ \rangle =:
\fW$ be the partial order defined as follows:
$$p\in{\fW}\ \Leftrightarrow\
  p\in\parto,$$
$$p\leq q \ \Leftrightarrow\ p\ceq^* q.$$

The distributivity number $\di (\fW)$
is defined as the least cardinal $\kappa$ for
which the Boolean algebra r.o.($\fW$) is not $\kappa$-distributive.
\hfill

\subsection*{The four cardinals are equal}

Now we will show, that the four cardinals defined above are all
equal. This is a similar result as in the case when we consider
infinite subsets of $\omega$ instead of infinite partitions
(cf.\,\cite{Plewik} and \cite{Balcar}).

\begin{fct}\label{fct:mao}
If $T\subseteq \parto$ is an open and dense set with respect to
the dual Ellentuck topology, then it contains a {\mao} family.
\end{fct}

\proof First choose an almost orthogonal family ${\cal A}\subseteq
T$ which is maximal in $T$. Now for an arbitrary $X\in\parto$,
$T\cap\partX\neq\emptyset$. So, $X$ must be compatible with
some $A\in {\cal A}$ and therefore ${\cal A}$ is {\mao}.
\eof

\begin{lm}\label{lm:h<add}
$\fH\leq\add (\fJ)$.
\end{lm}

\proof Let $\langle S_{\alpha}:\alpha<\lambda<\fH\rangle$ be a sequence
of nowhere dual-Ramsey sets and let $T_{\alpha}\subseteq\parto\setminus
S_{\alpha}$ ($\alpha<\lambda$) be such that $T_{\alpha}$ is open and dense
with respect to the dual Ellentuck topology (which is always possible by
the Remark 1). For each $\alpha<\lambda$
let $$T_{\alpha}^*:=\{X\in\parto:\exists Y\in T_{\alpha}(X\ceq^*
Y\wedge\neg(X\ass Y))\}.$$
It is easy to see, that for each $\alpha<\lambda$ the set
$T_{\alpha}^*$ is open and dense
with respect to the dual Ellentuck topology.\\
Let $U_{\alpha}\subseteq T_{\alpha}^*$ ($\alpha<\lambda$) be \mao.
Because $\lambda<\fH$, the set $\langle U_\alpha:\alpha<\lambda\rangle$
can not be shattering.
Let for $\alpha<\lambda$ $U_{\alpha}^*:=\{X\in\parto:\exists
Z_{\alpha}\in U_{\alpha}(X\ceq^* Z_{\alpha})\}$, then
$U_{\alpha}^*\subseteq T_{\alpha}$ and
$\bigcap_{\alpha<\lambda}U_{\alpha}^*$ is open and dense with
respect to the dual Ellentuck topology:\\
$\bigcap\limits_{\alpha<\lambda}U_{\alpha}^*$ is open: clear.\\
$\bigcap\limits_{\alpha<\lambda}U_{\alpha}^*$ is dense: Let
${\open{s}{Z}}$ be arbitrary. Because $\langle U_{\alpha}: \alpha
<\lambda\rangle$ is not shattering, there is a $Y\in {\open{s}{Z}}$ such
that $\forall \alpha <\lambda\exists X_{\alpha}\in
U_{\alpha}(Y\ceq^*X_{\alpha}).$ Hence,
$Y\in\bigcap_{\alpha<\lambda}U_{\alpha}^*$.
\hfill\smallskip

Further we have by construction
$$\bigcap\limits_{\alpha<\lambda}U_{\alpha}^*\cap
\bigcup\limits_{\alpha<\lambda}S_{\alpha} =\emptyset,$$
which completes the proof.
\eop

\begin{lm}\label{lm:h<d}
$\fH\leq \di (\fW)$.
\end{lm}

\proof Let $\langle T_\alpha :\alpha <\lambda <\fH\rangle$ be a
sequence of open and dense sets with respect to the dual Ellentuck
topology. Now the set $\bigcap_{\alpha <\lambda}U^*_\alpha$, constructed
as in Lemma~\ref{lm:h<add}, is dense (and even open) and a subset of
$\bigcap_{\alpha<\lambda}T_\alpha$. Therefore $\fH\leq\di (\fW)$.
\eop

\begin{lm}\label{lm:h>add}
$\add (\fJ)\leq\fH$.
\end{lm}

\proof Let $\langle R_{\alpha}:\alpha <\fH\rangle$ be a shattering family
and $P_{\alpha}:=\{X:\exists Y\in R_{\alpha}(X\ceq^*Y)\}$.\\
For each $\alpha <\fH$, $P_{\alpha}$ is dense and open with respect to
the dual Ellentuck topology:\\
$P_{\alpha}$ is open: clear.\\
$P_{\alpha}$ is dense: Let ${\open{s}{Z}}$ be arbitrary and $X\in
{\open{s}{Z}}$. Because $R_{\alpha}$ is \mao, there is a $Y\in
R_{\alpha}$ such that $X':=X\pint Y\in\parto$. Now let $X''\ass X'$ such
that $X''\in {\open{s}{Z}}$, then $X''\ceq^* Y$.\hfill\smallskip

Now we show that $\bigcap_{\alpha<\fH}P_{\alpha}=\emptyset$
and therefore
$\bigcup_{\alpha<\fH}(\parto\setminus P_\alpha)=\parto$.
Assume there is an
$X\in\bigcap_{\alpha<\fH}P_{\alpha}$, then
$\forall\alpha<\fH\exists Y_{\alpha}\in R_{\alpha}(X\ceq^* Y_{\alpha})$.
But this contradicts that $\langle R_{\alpha}:\alpha <\fH\rangle$ is
shattering.
\eop

\begin{lm}\label{lm:h>d}
$\di (\fW)\leq\fH$.
\end{lm}

\proof In the proof of Lemma~\ref{lm:h>add} we constructed a sequence
$\langle P_\alpha :\alpha<\fH\rangle$ of open and dense sets with
an empty intersection. Therefore $\bigcap_{\alpha<\fH}P_\alpha$ is not
dense.
\eop

\begin{cor}\label{cor:h>cov}
$\cov (\fJ)\leq\fH$.
\end{cor}

\proof In the proof of Lemma~\ref{lm:h>add}, in fact we proved that $\cov
(\fJ)\leq\fH$.
\eop

\begin{cor}\label{cor:a=c=h}
$\add (\fJ)=\cov (\fJ)=\di (\fW)=\fH$.
\end{cor}

\proof It is clear that $\add (\fJ)\leq \cov (\fJ)$. By the
Lemmas~\ref{lm:h<d} and \ref{lm:h>d} we know that $\fH=\di (\fW)$.
Further by the
Lemma~\ref{lm:h<add} and the Corollary~\ref{cor:h>cov} it follows that
$\fH\leq\add (\fJ)\leq\cov (\fJ)\leq\fH$. Hence we have $\add (\fJ)=\cov
(\fJ)=\di (\fW)=\fH$.
\eop

\begin{cor}
The union of less than $\fH$ completely dual-Ramsey sets is dual-Ramsey,
but the union of $\fH$ completely dual-Ramsey sets can be a set, which
does not have the dual-Ramsey property.
\end{cor}

\proof
Follows from Remark~1 and Corollary~\ref{cor:a=c=h}.
\eop

\subsection*{On the consistency of $\fH \bold{>\omega_1}$}

First we give some facts concerning the dual-Mathias forcing.
\hfill\smallskip

The conditions of dual-Mathias forcing are pairs $\langle s,X\rangle$
such that $s\in\NN$, $X\in\parto$ and $s\ceq X$,
stipulating $\langle s,X\rangle
\leq\langle t,Y\rangle$ if and only if ${\open{s}{X}}\subseteq {\open{t}
{Y}}$.\\
It is not hard to see that similar to Mathias forcing,
the dual-Mathias forcing can
be decomposed as $\fW * \Padname$, where $\fW$ is defined as above and
$\Padname$ denotes dual-Mathias forcing with conditions only with second
coordinate in $\tilde{\fD}$,
where $\tilde{\fD}$ is an $\fW$-generic object.\\
Further, because dual-Mathias forcing has pure decision (cf.\,\cite{dual}),
it is proper and has the Laver property and therefore adds no Cohen reals.
\hfill\medskip

If we make an $\omega_2$-iteration of dual-Mathias forcing with countable
support, starting from a model in which the continuum hypothesis holds,
we get a model in which the dual-shattering cardinal $\fH$ is equal to
$\omega_2$.
\hfill\smallskip

Let $V$ be a model of CH and let
$P_{\omega_2}:= \langle P_\alpha,\dot{Q}_\beta :
\alpha\leq\omega_2,\beta <\omega_2\rangle$ be a countable support
iteration of dual-Mathias forcing, i.e.\,$\forall\alpha<\omega_2:
\force_{P_\alpha}``\dot{Q}_\alpha\,\mbox{is dual-Mathias forcing}"$.
\hfill\smallskip

In the sequel we will not distinguish between a member of $\fW$ and
its representative. In the proof of the following theorem, a set
$C\subseteq\omega_2$ is called $\omega_1$-club if $C$ is unbounded in
$\omega_2$ and closed under increasing sequences of length $\omega_1$.

\begin{thm}\label{thm:h>o1}
If $G$ is $P_{\omega_2}$-generic over $V$, where $V\models {\operatorname
{CH}}$, then $V[G]\models\fH=\omega_2$.
\end{thm}

\proof In $V[G]$ let $\langle D_\nu :\nu <\omega_1\rangle$ be a family of
open dense subsets of $\fW$. Because dual-Mathias forcing is proper
and by a standard L\"owenheim-Skolem argument, we find
a $\omega_1$-club $C\subseteq\omega_2$ such that for each $\alpha\in C$
and every $\nu <\omega_1$ the set $D_\nu\cap V[G_\alpha]$ belongs to
$V[G_\alpha]$ and is open dense in $\fW^{V[G_\alpha]}$. Let $A\in\fW^{V
[G]}$ be arbitrary. By properness and genericity and because
$P_{\omega_2}$ has countable support, we may assume that $A\in G(\alpha)'$
for an $\alpha\in C$, where $G(\alpha)'$ is the first component according
to the decomposition of Mathias forcing
of the $\dot{Q}_\alpha [G_\alpha]$-generic object
determined by $G$. As $\alpha\in C$, $G(\alpha)'$ clearly meets every
$D_\nu$ ($\nu <\omega_1$). But now $X_\alpha$, the $\dot{Q}_\alpha$-generic
partition (determined by $G(\alpha)''$) is below each member of $G(\alpha)'$,
hence below $A$ and in $\bigcap_{\nu<\omega_1}D_\nu$. Because $A$ was
arbitrary, this proves that $\bigcap_{\nu<\omega_1}D_\nu$ is dense in
$\fW$ and therefore $\di (\fW)>\omega_1$. Again by properness of dual-Mathias
forcing $V[G]\models 2^{\omega_0}=\omega_2$ and we finally have
$V[G]\models \fH=\omega_2$.
\eop

In the model constructed in the proof of Theorem~\ref{thm:h>o1} we
have $\fH >\ft$, where $\ft$ is the well-known tower number (for a
definition of $\ft$ cf.\,\cite{vDouwen}). Moreover, we can show

\begin{cor} The statement $\fH >\cov ({\cal M})$ is relatively consistent
with $\operatorname{ZFC}$, (where ${\cal M}$ denotes the ideal of meager
sets).
\end{cor}

\proof Because dual-Mathias forcing is proper and does not add Cohen reals,
also forcing with $P_{\omega_2}$ does not add Cohen reals. Further it is
known that $\ft\leq \cov ({\cal M})$ (cf.\,\cite{Piotro} or \cite{book}).
Now because forcing with $P_{\omega_2}$
does not add Cohen reals, in $V[G]$ the covering number $\cov(
{\cal M})$ is still $\omega_1$ (because each real in $V[G]$ is in a
meager set with code in $V$). This completes the proof.
\eop

\rmktwo In \cite{vDouwen} Theorem\,3.1.(c) it is shown that
$\omega\leq\kappa <\ft$ implies that $2^{\kappa}=2^{\omega_0}$.
We do not have a similar result for the dual-shattering cardinal
$\fH$. If we start our forcing construction $P_{\omega_2}$
with a model $V\models \operatorname{CH} + 2^{\omega_1}
=\omega_3$, then (again by properness of dual-Mathias forcing)
$V[G]\models \fH =\omega_2 =2^{\omega_0} <2^{\omega_1}=\omega_3$
(where $G$ is $P_{\omega_2}$-generic over $V$).

{\bf Remark:} Recently Spinas showed in \cite{spinas}, that
$\fH<\fh$ is consistent with ZFC. But it is still open
if MA$+\neg$CH implies that $\omega_1 <\fH$.

\section{On the dual-splitting cardinals $\fS$ and $\fS'$}\label{sec:s}

Let $X_1,X_2$ be two partitions. We say $X_1$ {\em splits\/} $X_2$ if
$X_1\comp X_2$ and it exists a partition $Y\ceq X_2$, such that
$X_1\ort Y$. A family ${\cal S}\subseteq\parto$ is called
{\em splitting\/} if for each non-trivial $X\in\parto$ there exists
an $S\in{\cal S}$ such that $S$ splits $X$. The dual-splitting cardinal $\fS$
(resp. $\fS'$) is the least  cardinal number $\kappa$, for which there exists
a splitting family ${\cal S}\subseteq\parto$ (resp. ${\cal S}
\subseteq\partou$) of cardinality $\kappa$.\hfill\smallskip

It is obvious that $\fS\leq\fS'$.\hfill\bigskip

First we compare the dual-splitting number $\fS'$ with
the well-known bounding number $\fb$ (a definition of $\fb$ can
be found in \cite{vDouwen}).

\begin{thm}\label{thm:bleqS}
$\fb\leq\fS'$.
\end{thm}

\proof Assume there exists a family ${\cal S}=\{S_\iota:\iota <\kappa
<\fb\}\subseteq\partou$ which is splitting. Let $B=\{b_\iota:\iota<
\kappa\}\subseteq\reals$ a set of infinite subsets of $\omega$
such that $b_\iota\in S_\iota$ (for all $\iota <\kappa$). Let
$f_{b_\iota}\in\oto$ be the (unique) increasing function such that
range($f_{b_\iota}$)=$b_\iota$. Because $\kappa<\fb$, the set
$\{f_{b_\iota}:\iota<\kappa\}$ is not unbounded. Therefore there
exists a function $d\in\oto$ such that $f_{b_\iota}<^*d$ (for all
$\iota<\kappa$). Now with the function $d$ we construct an
infinite partition $D$. First we define an infinite set of pairwise
disjoint finite sets $p_i$ ($i\in\omega$):
$$p_i:=[d^i(0),d^{i+1})$$
where $d^i$ denote the $i$-fold composition of $d$.\hfill\smallskip

Now the blocks of $D$ are defined as follows:
$$n\ \mbox{is in the $k$th block of $D$ {\em iff}}\ n\in p_i\wedge
i-\mmax\{\frac{l}{2}(l+1)<i:l\in\omega\}=k.$$
Because $d$ dominates $B$, for all $b_\iota\in B$ there exists a
natural number $m_\iota$, such that for all $i>m_\iota$:
$d^i(0)\leq b_\iota(d^i(0))<d^{i+1}(i)$ (cf.\,\cite{vDouwen}\,p.\,121).
So, for all $i>m_\iota$, $p_i\cap b_\iota\neq\emptyset$ and
therefore by the construction of the blocks of $D$, $b_\iota$ intersects
each block of $D$. But this implies, that $D$ is not compatible with
any element of ${\cal S}$ and so ${\cal S}$ can not be a splitting
family.
\eop

\begin{cor}
It is consistent with $\operatorname{ZFC}$, that $\fs <\fS'$.
\end{cor}

\proof Because $\fb\leq\fS'$ is provable in ZFC, it is enough to
prove that $\fs<\fb$ is consistent with ZFC, which is proved
in \cite{Shelah}.
\eop

Now we show that $\cov ({\cal M})\leq\fS$ (where ${\cal M}$
denotes the ideal of meager sets).
In \cite{dia} it is shown that if $\kappa<\cov ({\cal M})$ and
$\{X_\alpha:\alpha <\kappa\}\subseteq\parto$ is a family of
partitions, then there exists $Y\in\parto$ such that $Y\ort X_\alpha$
for each $\alpha <\kappa$. This implies the following

\begin{cor}\label{cor:covMleqS}
$\cov ({\cal M})\leq \fS$.
\end{cor}

\proof Let $S,Y\in\parto$. If $S\ort Y$, then $S$ does not split $Y$
and therefore a family of cardinality less than $\cov ({\cal M})$
can not be splitting.
\eop

As a corollary we get again a consistency result.

\begin{cor}
It is consistent with $\operatorname{ZFC}$, that $\fs <\fS$.
\end{cor}

\proof If we make an $\omega_1$-iteration of Cohen forcing with
finite support starting from a model $V\models \cov
({\cal M})=\omega_2=\fc$, we get a model in which $\omega_1=\fs
<\cov ({\cal M})=\omega_2=\fc$ holds. Hence, by
Corollary~\ref{cor:covMleqS}, this is a model for $\omega_1=\fs
<\fS=\omega_2$.
\eop

Until now we have $\cov ({\cal M}),\fb\leq\fS'$, which would be
trivial if one could show that $\fS'=\fc$, where $\fc$ is the
cardinality of ${\cal P}(\omega)$. But this is not the case
(cf.~\cite{dia}). For the sake of completeness we will give now
the notion of forcing used in \cite{dia} to construct
a model in which we have $\fS' <\fc$.
\hfill\smallskip

Let $\oQ$ the notion of forcing defined as follows. The conditions
of $\oQ$ are pairs $\langle s,A\rangle$ such that $s\in\NN$, $A\in
\partf$ and $s\seg A$, stipulating $\langle s,A\rangle \leq \langle
t,B\rangle$ if and only if $t\seg s$ and $B\ceq A$.
($s$ is called the stem of the condition.) If $\langle
s,A_1\rangle ,\langle s,A_2\rangle$ are two $\oQ$-conditions, then
$\langle s,A_1\pint A_2\rangle\leq\langle s,A_1\rangle ,\langle
s,A_2\rangle$. Hence, two $\oQ$-conditions with the same stem are
compatible and because there are only countably many stems, the
forcing notion $\oQ$ is $\sigma$-centered.\\
Now we will see, that forcing with $\oQ$ adds an infinite partition
which is compatible with all old infinite partitions but is not
contained in any old partition. (So, the forcing notion $\oQ$ is in
a sense like the dualization of Cohen forcing.)

\begin{lm}\label{lm:Q}
If $G$ is $\oQ$-generic over $V$, then $G\in\partou$ and
$V[G]\models\forall X\in\parto\cap V (G\comp X\wedge \neg (X\ceq^*
G)).$
\end{lm}

\proof Let $X\in V$ be an arbitrary, infinite partition. The set
$D_n$ of $\oQ$-conditions $\langle s,A\rangle$, such that
\begin{itemize}
\item[  (i)] at least one block of $s$ has more than $n$ elements,
\item[ (ii)] at least $n$ blocks of $X$ are each the union of blocks
of $A$,
\item[(iii)] there are at least $n$ different blocks $b_i\in X$, such
that $\bigcup b_i\in s\puni X$,
\end{itemize}
is dense in $\oQ$ for each $n\in\omega$. Therefore, at least one
block of $G$ is infinite (because of (i)), $G$ is compatible
with $X$ (because of (ii)) and $X$ is not coarser${}^*$ than $G$
(because of (iii)). Now, because $X$ was arbitrary, the $\oQ$-generic
partition $G$ has the desired properties.
\eop

Because the forcing notion $\oQ$ is $\sigma$-centered and each
$\oQ$-condition can be encoded by a real number, forcing with
$\oQ$ does neither collapse any cardinals nor change the cardinality
of the continuum and we can prove the following

\begin{lm}
It is consistent with $\operatorname{ZFC}$ that $\fS'<\fc$.
\end{lm}

\proof \cite{dia} If make an $\omega_1$-iteration of $\oQ$ with
finite support, starting from a model in which we have $\fc=\omega_2$,
then the $\omega_1$ generic objects form a splitting family.
\eop
{}\vspace{0.5cm}

Even if a partition does not have a complement, for each non-trivial
partition $X$ we can define a non-trivial partition $Y$, such
that $X\ort Y$.\hfill\smallskip

Let $X=\{b_i:i\in\omega\}\in\parto$ and assume
that the blocks $b_i$ are ordered by their least element and
that each block is ordered by the natural order. A block
is called trivial, if it is a singleton. With respect
to this ordering define for each non-trivial partition $X$ the
partition $X^{\Lort}$ as follows.
\hfill\smallskip

If $X\in\partou$ then
\begin{center}
$n$ is in the $i$th block of $X^{\Lort}$\\
{\em iff}\\
$n$ is the $i$th element of a block of $X$,
\end{center}
otherwise
\begin{center}
$n,m$ are in the same block of $X^{\Lort}$\\
{\em iff}\\
$n,m$ are both least elements of blocks of $X$.
\end{center}

It is not hard to see that for each non-trivial $X\in\parto$,
$X\ort X^\Lort$.\hfill\smallskip

A family ${\cal W}\subseteq\partou$ is called {\em weak splitting},
if for each partition $X\in\parto$, there is a $W\in{\cal W}$ such
that $W$ splits $X$ or $W$ splits $X^\Lort$.
The cardinal number $w\fS$ is the least cardinal
number $\kappa$, for which there exists a
weak splitting family of cardinality $\kappa$.
(It is obvious that $w\fS\leq\fS'$.)\\
A family $U$ is called a $\pi$-{\em base\/} for a free
ultra-filter ${\cal F}$ over $\omega$ provided for every
$x\in{\cal F}$ there
exists $u\in U$ such that $u\subseteq x$. Define
$$\pi\fu :=\mmin \{|U|:U\subseteq\reals\,\mbox{is a
$\pi$-base for a free ultra-filter over $\omega$}\}.$$
In \cite{BalcarSimon} it is proved, that $\pi\fu =\fr$ (see also
\cite{Vaughan} for more results concerning $\fr$).
\hfill\smallskip

Now we can give an upper and a lower bound for the size of $w\fS$.

\begin{thm}\label{thm:wSlequ}
$w\fS\leq\fr$.
\end{thm}

\proof We will show that $w\fS\leq\pi\fu$.
Let $U:=\{u_\iota\in\reals: \iota<\pi\fu\}$ be a $\pi$-basis for a
free ultra-filter ${\cal F}$ over $\omega$.
W.l.o.g.\quad we may assume,
that all the $u_\iota\in U$ are co-infinite.
Let ${\cal U}=\{Y_u\in\parto : u\in U \wedge Y_u=\{u_i: u_i=u\vee
(u_i=\{n\}\wedge n\not\in u)\}\}$. Now we take an arbitrary $X=
\{b_i: i\in\omega\}\in\parto$ and define for every $u\in U$
the sets $I_u:=\{i:b_i\cap u\neq\emptyset\}$ and $J_u:=\{j:b_j\cap
u=\emptyset\}$. It is clear that $I_u\cup J_u=\omega$ for every $u$.
\hfill\smallskip

If we find a $u\in U$ such that $|I_u|=|J_u|=\omega$, then $Y_u$
splits $X$. To see this, define the two infinite partitions
$$Z_1:=\{a_k: a_k=\bigcup_{i\in I_u}b_i\,\vee\,\exists j\in J_u
a_k=b_j\}$$
and
$$Z_2:=\{a_k: a_k=\bigcup_{j\in J_u}b_j\,\vee\,\exists i\in I_u
a_k=b_i\}.$$
Now we have $X\puni Y_u=Z_1$ (therefore $Z_1\ceq X,Y_u$) and $Z_2
\ceq X$ but $Z_2\ort Y_u$.\\
(If each block of $b_i$ is finite, then we are always in this case.)
\hfill\smallskip

If we find an $x\in{\cal F}$ such that $|I_x|<\omega$ (and therefore
$|J_x|=\omega$), then we find an $x'\subseteq x$, such that $|I_x|=1$
and for this $i\in I_x$, $|b_i\setminus x'|=\omega$. (This is because
${\cal F}$ is a free ultra-filter.) Now take a $u\in U$ such that
$u\subseteq x'$ and we are in the former case for $X^\Lort$.
Therefore, $Y_u$ splits $X^\Lort$.\hfill\smallskip

If we find an $x\in{\cal F}$ such that $|J_x|<\omega$ (and therefore
$|I_x|=\omega$), let $I(n)$ be an enumeration of $I_x$ and define
$y:=x\cap\bigcup_{k\in\omega}b_{I(2k)}$. Then $y\subseteq x$ and
$|x\setminus y|=\omega$. Hence, either $y$ or $\omega\setminus y$
is a superset of some $u\in U$. But now $|J_u|=\omega$ and
we are in a former case.
\eop

A lower bound for $w\fS$ is $\cov ({\cal M})$.

\begin{thm}
$\cov ({\cal M})\leq w\fS$.
\end{thm}

\proof Let $\kappa <\cov ({\cal M})$ and ${\cal W}=\{W_\iota: \iota
<\kappa\}\subseteq\partou$. Assume for each $W_\iota\in {\cal W}$
the blocks are ordered by their least element and each block is
ordered by the  natural order. Further assume that $b_{i(\iota )}$
is the first block of $W_\iota$ which is infinite.
Now for each $\iota <\kappa$ the set
$D_\iota$ of functions $f\in\oto$ such that

$$\begin{array}{ll}
\forall n,m,k\in\omega &\exists h\in\omega
t_1\in b_n,t_2\in b_m,t_3,t_4\in b_h\exists s\in b_{i(\iota)}\\
     &f(t_1)=f(t_3)\wedge f(t_2)=f(t_4)\wedge
      |\{s'\leq s: f(s')=f(s)\}|=k+1.
\end{array}$$

is the intersection of countably many open dense sets and therefore
the complement of a meager set. Because $\kappa <\cov ({\cal M})$,
we find an unbounded function $g\in\oto$ such that $g\in\bigcap_{\iota
<\kappa}D_\iota$. The partition $G=\{g^{-1}(n): n\in\omega\}\in\partou$
is orthogonal with each member of ${\cal W}$ and for each $W_\iota
\in {\cal W}$ and each $k\in\omega$, there exists an $s\in b_{i(\iota)}$,
such that $s$ is the $k$th element of a block of $G$.
Hence, ${\cal W}$ can not be a weak splitting family.
\eop

\section{On the dual-reaping cardinals $\fR$ and $\fR'$}

A family ${\cal R}\subseteq\parto$ is called {\em reaping} (resp.
{\em reaping}\,$'$), if
for each partition $X\in\parto$ (resp. $X\in\partou$) there exists a
partition $R\in {\cal R}$ such that $R\ort X$ or $R\ceq^* X$.
The dual-reaping cardinal $\fR$ (resp. $\fR'$) is the least cardinal
number $\kappa$, for which there exists a
reaping (resp. reaping$'$) family of cardinality $\kappa$.
\hfill\smallskip

It is clear that $\fR'\leq\fR$. Further by finite modifications of
the elements of a reaping family, we may replace $\ceq^*$ by $\ceq$
in the definition above.\hfill\bigskip

If we cancel in the definition of the reaping number the expression
$``R\ceq^* X"$, we get the definition of an orthogonal family.
\hfill\smallskip

A family ${\cal O}\subseteq\parto$ is called {\em orthogonal} (resp.
{\em orthogonal}\,$'$), if for each non-trivial partition $X\in\parto$
(resp. for each partition $X\in\partou$)
there exists a partition $O\in {\cal O}$ such that $O\ort X$. The
dual-orthogonal cardinal $\fO$ (resp. $\fO'$) is the least  cardinal
number $\kappa$, for which there exists a  orthogonal (resp.
orthogonal$'$) family of cardinality $\kappa$. (It is obvious
that $\fO'\leq\fO$.) Note, that $\fo =\fc$, where $\fc$ is the
cardinality of ${\cal P}(\omega)$ and $\fo$ is defined
like $\fO$ but for infinite subsets of $\omega$ instead of infinite
partitions. (Take the complements of a maximal antichain in $\reals$
of cardinality $\fc$. Because an orthogonal family must avoid all
this complements, it has at least the cardinality of this maximal
antichain.)\hfill\smallskip

It is also clear that each orthogonal${}^( {'}{}^)$ family is also
a reaping${}^( {'}{}^)$ family and therefore
$\fR^( {'}{}^) \leq\fO^( {'}{}^)$. Further one can show that $\fR'$ is
uncountable (cf.~\cite{dia}).
Now we show that $\fO'\leq\fd$, where $\fd$
is the well-known dominating number
(for a definition cf.~\cite{vDouwen}), and that
$\cov ({\cal M})\leq\fO'$.

\begin{lm}
$\fO'\leq\fd$.
\end{lm}

\proof Let $\{d_\iota\in\oto: \iota<\fd\}$ be a dominating family.
Then it is not hard to see that the family $\{D_\iota:\iota<\kappa\}
\subseteq\parto$, where each $D_\iota$ is constructed from $d_\iota$
like $D$ from $d$ in the proof of Theorem~\ref{thm:bleqS}, is an
orthogonal family.
\eop

Let $\fI$ be the least cardinality of an independent family (a
definition and some results can be found in \cite{Kunen}), then

\begin{lm}
$\fO\leq\fI$.
\end{lm}

\proof Let $I\subseteq\reals$ be an independent family of cardinality
$\fI$. Let $I':=\{r\in\reals :r\ass\bigcap {\cal A}\setminus\bigcup
{\cal B}\}$, where ${\cal A},{\cal B}\in [I]^{<\omega}$, ${\cal A}\neq
\emptyset$, ${\cal A}\cap {\cal B}=\emptyset$ and $r\ass x$ means
$|(r\setminus x)\cup (x\setminus r)|<\omega$. It is not hard to see
that $|I'|=|I|=\fI$. Now let ${\cal I}={\cal I}_1\cup {\cal I}_2$
where ${\cal I}_1:=\{X_r\in\parto :r\in I'\wedge X_r=\{b_i: b_i=r\vee
(b_i=\{n\}\wedge n\not\in r)\}\}$ and ${\cal I}_2:=\{Y_r:\exists X_r
\in {\cal I}_1 (Y_r=X_r^\Lort)\}$. We see, that ${\cal I}\subseteq
\parto$ and $|{\cal I}|=\fI$. It leave to show that ${\cal I}$ is an
orthogonal family.\\
Let $Z\in\parto$ be arbitrary and let $r:=\mMin (Z)$.
If $r\in I'$, then $X_r\ort Z$ (where $X_r\in {\cal I}_1$). And if
$r\not\in I'$, then there exists an $r'\in I'$ such that $r\cap r'=
\emptyset$. But then $Y_{r'}\ort Z$ (where $Y_{r'}\in {\cal I}_2$).
\eop

Because $\fR\leq\fO$, the cardinal number $\fI$ is also an upper
bound for $\fR$. But for $\fR$, we also find another upper bound.

\begin{lm}
$\fR\leq\fr$.
\end{lm}

\proof Like in Theorem~\ref{thm:wSlequ} we show that $\fR\leq\pi\fu$.
Let $U:=\{u_\iota\in\reals: \iota<\pi\fu\}$ be as in the proof
of Theorem~\ref{thm:wSlequ} and let
${\cal U}=\{Y_u\in\parto : u\in U \wedge Y_u=\{u_i: u_i=\omega
\setminus u\vee (u_i=\{n\}\wedge n\in u)\}\}$. Take an arbitrary
partition $X\in\parto$. Let $r:=\mMin (X)$ and $r_1:=\{n\in r:
\{n\}\in X\}$. If we find a $u\in U$ such that $u\subseteq r_1$,
then $Y_u\ceq X$. Otherwise, we find a $u\in U$ such that either
$u\subseteq\omega\setminus r$ or $u\subseteq r\setminus r_1$ and
in both cases $Y_u\ort X$.
\eop

Now we will show, that it is consistent with $\operatorname{ZFC}$
that $\fO$ can be small. For this we first show, that a Cohen real
encode an infinite partition which is orthogonal to each old
non-trivial infinite partition. (This result is in fact a
corollary of Lemma~5 of \cite{dia}.)

\begin{lm}\label{lm:C}
If $c\in\oto$ is a Cohen real over $V$, then $C:=\{c^{-1}(n):n\in
\omega\}\in\partou\cap V[c]$ and $\forall X\in\parto\cap V
(\neg (X\ass\{\omega \})\rightarrow C\ort X)$.
\end{lm}

\proof We will consider the Cohen-conditions as finite sequences
of natural numbers, $s=\{s(i): i<n<\omega\}$.
Let $X=\{b_i: i\in\omega\}\in V$ be an arbitrary, non-trivial infinite
partition. The set $D_{n,m}$ of Cohen-conditions $s$, such that
\begin{itemize}
\item[  (i)] $|\{i: s(i)=0\}|\geq n$,
\item[ (ii)] $\exists k>n \exists i(s(i)=k)$,
\item[(iii)] $\exists a_n\in b_n\exists a_m\in b_m\exists l\exists
a_1,a_2\in b_l(s(a_n)=s(a_1)\wedge s(a_m)=s(a_2))$,
\end{itemize}
is a dense set for each $n,m\in\omega$. Now, because $X$ was
arbitrary, the infinite partition $C$ is orthogonal to each infinite
partition which is in $V$. (Note that because of (i), $C\in\partou$.)
\eop

Now we can show, that $\fO$ can be small.

\begin{lm}\label{lm:smallO}
It is consistent with $\operatorname{ZFC}$ that $\fO<\cov ({\cal M})$.
\end{lm}

\proof If make an $\omega_1$-iteration of Cohen forcing with
finite support, starting from a model in which we have
$\fc=\omega_2=\cov ({\cal M})$,
then the $\omega_1$ generic objects form an orthogonal family.
Now because this $\omega_1$-iteration of Cohen forcing does not
change the cardinality of $\cov ({\cal M})$, we have a model in
$\omega_1 =\fO<\cov ({\cal M})=\omega_2$ holds.
\eop

Because $\fR\leq\fO$ we also get the relative consistency
of $\fR<\cov ({\cal M})$. Note that this is not true for $\fr$.

As a lower bound for $\fR'$ we find $\fp$, where $\fp$ is the
pseudo-intersection number (a definition of $\fp$ can be found
in \cite{vDouwen}).

\begin{lm}
$\fp\leq\fR'$.
\end{lm}

\proof In \cite{Bell} it is proved that $\fp=\frak{m}_{\sigma
\mbox{\footnotesize{-centered}}}$, where $$\frak{m}_{\sigma
\mbox{\footnotesize{-centered}}}=\mmin \{\kappa :``\MA
(\kappa )\mbox{\,for $\sigma$-centered posets" fails\,}\}.$$
Let ${\cal R}=\{ R_\iota :\iota <\kappa <\fp\}$
be a set of infinite partitions.
Now remember that the forcing notion $\oQ$ (defined in
section~\ref{sec:s}) is $\sigma$-centered and because $\kappa <\fp$
we find an $X\in\partou$ such that ${\cal R}$ does not reap $X$.
\eop

As a corollary we get

\begin{cor}
If we assume $\MA$, then $\fR'=\fc$.
\end{cor}

\proof If we assume $\MA$, then $\fp=\fc$.
\eop

\section{What's about towers and maximal (almost) orthogonal families?}

Let $\kappa_{mao}$ be the least cardinal
number $\kappa$, for which there
exists an infinite {\mao} family of cardinality $\kappa$. And let
$\kappa_{tower}$ be the least cardinal number $\kappa$,
for which there exists
a family ${\cal F}\subseteq\parto$ of cardinality $\kappa$, such
that ${\cal F}$ is well-ordered by $\ceq^*$ and $\neg\exists Y\in
\parto\forall X\in {\cal F}(Y\ceq^* X)$.\hfill\smallskip

Now Krawczyk proved that $\kappa_{mao}=\fc$ (cf.\,\cite{dia}) and
Carlson proved that $\kappa_{tower}=\omega_1$ (cf.\,\cite{Matet}).
So, these cardinals do not look interesting.
But what happens if we cancel the word
``almost" in the definition of $\kappa_{mao}$?\\
A family ${\cal F}\subseteq
\parto$ (resp. ${\cal F}\subseteq\partou$) is a {\em maximal
anti-chain\/} in $\parto$ (resp. $\partou$), if ${\cal F}$ is a maximal
infinite family of pairwise orthogonal partitions. Let $\kappa_A$ (resp.
$\kappa_{A'}$) be the least cardinality of a maximal anti-chain in
$\parto$ (resp. $\partou$).\\
Note that the corresponding cardinal for infinite
subsets of $\omega$ would be equal to $\omega$.\hfill\smallskip

First we know that $\cov ({\cal M})\leq \kappa_A,\kappa_{A'}$
(which is proved in
\cite{dia}) and $\fb\leq \kappa_{A'}$ (which one can prove like
Theorem~\ref{thm:bleqS}). Further it is not hard to see that
$\kappa_A\leq \kappa_{A'}$.\\
But these results concerning $\kappa_A$ and $\kappa_{A'}$
are also not interesting,
because Spinas showed in \cite{spinas} that
$\kappa_A = \kappa_{A'} = \fc$.


\section{The diagram of the results}

Now we summarize the results proved in this article
together with other known results.
\hfill\medskip

{\em splitting:}
$$\begin{array}{ccccccc}
\fb & \rule[0.8ex]{10ex}{0.1ex} & \fS' &\hspace{-5ex}\rule[0.8ex]{25ex}{0.1ex}
\hspace{-20ex} & & & \fc \\
\rule[0.5ex]{0.1ex}{2ex}& &\rule[-4.5ex]{0.1ex}{7ex} & & &
&\rule[-10.5ex]{0.1ex}{13ex}
\vspace{-10.5ex}

\\
\fh & &     &                &      &                &     \\
\rule[0ex]{0.1ex}{2ex}    &  &     &    &      &  &     \\
\fH & \rule[0.8ex]{10ex}{0.1ex} & \fS & &      &                &     \\
\rule[0ex]{0.1ex}{7ex}& & \rule[0ex]{0.1ex}{7ex}& & & & \rule{0.1ex}{6ex}  \\
\omega_1 & \rule[0.8ex]{10ex}{0.1ex} & \cov ({\cal M}) &
\rule[0.8ex]{10ex}{0.1ex} & w\fS
& \rule[0.8ex]{10ex}{0.1ex} & \fr
\end{array}$$

{\samepage
{\em reaping:}
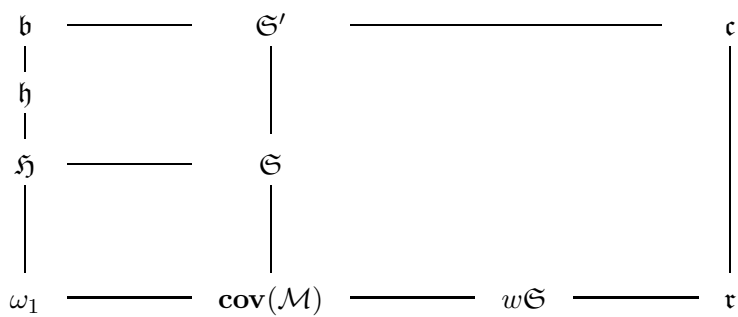
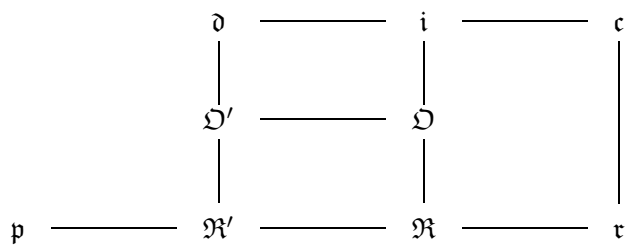
$$\begin{array}{ccccccc}
{ }& & \fd & \rule[0.8ex]{10ex}{0.1ex} & \fI &\rule[0.8ex]{10ex}{0.1ex}&\fc \\
  &  &\rule[0ex]{0.1ex}{5ex} & &\rule[0ex]{0.1ex}{5ex} & &
\rule[-8ex]{0.1ex}{13ex}\vspace{-8ex}

\\
    & & \fO' & \rule[0.8ex]{10ex}{0.1ex} & \fO &           &     \\
  &  &\rule[0ex]{0.1ex}{5ex} & &\rule[0ex]{0.1ex}{5ex} & & \\
\fp & \rule[0.8ex]{10ex}{0.1ex} &\fR' &\rule[0.8ex]{10ex}{0.1ex}& \fR
&\rule[0.8ex]{10ex}{0.1ex} & \fr
\end{array}$$
(In the diagrams, the invariants grow larger, as one moves up or to the
right.)}

\subsubsection*{Consistency results:}
\begin{itemize}
\item $\cov ({\cal M}) <\fH$\ ; $\fH < \fh$\ ;
      \ $\fH <\cov ({\cal M})$ (this is
      because $\fh <\cov ({\cal M})$ is consistent with ZFC)
\item $\fs<\fS$\ ;\ $\fS'<\fc$
\item $\fO<\cov ({\cal M})$
\end{itemize}

{\sc Note added in proof:} Recently, J\"org Brendle informed me that
he has proved, that $\MA\,+\,\fH\,<\,2^{\aleph_0}$ is consistent with ZFC.

\medskip

Lorenz Halbeisen\\
Departement Mathematik\\
ETH-Zentrum\\
8092 Z\"urich\\
Switzerland\hfill\smallskip

E-mail: halbeis@@math.ethz.ch

\end{document}